\input amstex

\documentstyle{amsppt}
\topmatter
\title
Solitons and projectively flat affine surfaces
\endtitle
\author
Wlodzimierz Jelonek
\endauthor
\abstract{The aim of this paper is to give a  local
description of affine surfaces, whose induced Blaschke structure is
projectively flat. We show that such affine surfaces with constant Gauss
affine curvature and indefinite induced Blaschke metric are described
by soliton equations.}
\endabstract
\endtopmatter
\document

\define\hn{\hat \nabla}
\define\nh{\hat\nabla}
\define\n{\nabla}
\define\str{(\n,h,S)}

\define\p{\partial}
\define\e{\epsilon}
\define\lb{\lambda}
\define\a{\alpha}
\par
\medskip
{\bf 0. Introduction.} It is well known that the sine-Gordon (S-G) equation
$\p_{12}u=\text{sin} u $ is strictly related with Riemannian surfaces of
constant negative Gauss curvature immersed the Euclidean space
$\Bbb{R}^3$. The equation S-G has ben studied in this context by many
geometers since the end of the nineteenth century. This hyperbolic nonlinear
equation has soliton solutions. All the surfaces of constant negative
curvature immersed in $\Bbb{R}^3$ are not convex. We show in the present
paper that the similar soliton equations appear in the natural way in the
study of non-convex projectively flat affine surfaces with constant
affine Gauss curvature.  Among them are semi-Riemannian surfaces of
constant Gauss curvature immersed in the Lorentz space $\Bbb{R}^3$.
All the
soliton equations obtained for non-convex surfaces are of great
importance in the soliton theory (see [B-C]).
In the first part of the paper we consider affine locally symmetric
surfaces (with $\Bbb{C}$-diagonalizable shape operator). The similar
results can be found in ([S], [P]). The second part is devoted to
description of affine spheres (see [J-2]). The similar results were
obtained by V.V.Nesterenko [N]  and U.Simon and C.Wang in [S-W].

\par
\medskip
{\bf 1.  Preliminaries. } Our notation is as in \cite{N-P-2},
\cite{J-1}. Let $(M,f)$ be an affine surface in $\Bbb{R}^3$
with the induced equiaffine structure $\str$. It means that
$f:M\rightarrow \Bbb{R}^3$ is an immersion and there exists
a transversal section
$\xi$ of the vector bundle $f^*T\Bbb{R}^3$ such that
$$D_Xf_*Y=f_*(\n_XY)+h(X,Y)\xi$$
and $D_X\xi=-f_*(SX)$. We call $h,S$ the induced affine metric and
shape operator respectively. We shall assume that an immersion
$f$ is nondegenerate (that is $h$ is a nondegenerate semi-metric).
By $C$ we denote the cubic
form $C=\n h$.
By $J =
\frac18 h(C,C)$ we denote the Fubini-Pick invariant.  In
any basis in $TM$ we have
 $ h(C,C) = \sum  h^{ip}h^{jq}h^{kr}C_{ijk}C_{pqr}$.
For an equiaffine structure $(\nabla ,h,S)$ the following equations hold:
$$\gather
 R(X,Y)Z = h(Y,Z)SX - h(X,Z)SY \tag Gauss
\\ C(X,Y,Z) = C(Y,X,Z) \tag Codazzi 1
\\ \nabla S(X,Y) = \nabla S(Y,X) \tag Codazzi 2
\\ \hn S(X,Y) - \hn S(Y,X) = K(SX,Y)-K(SY,X) \tag C2'
\\ h(SX,Y) = h(Y,SX), \tag Ricci
\endgather$$
where $\hn$ is the Levi-Civita connection of the affine metric $h$
and $K$, defined by $K(X,Y)=\n_XY-\hn_XY$, is the difference tensor.
The Radon theorem says that the structures $\str$ on surface $M$
satisfying
the above equations  are in one-to-one
correspondence (up to an affine transformation)
with equiaffine surfaces $(M,f,\xi)$ in $\Bbb{R}^3$.
Equation (Codazzi 1) is equivalent to
$$h(K_XY,Z)=h(Y,K_XZ).  \tag R1$$
The Blaschke structure $\str$ is characterized additionally by the
apolarity condition
$$\text{tr}_hK=0.  \tag A$$
For a Blaschke structure on a surface $M$ we denote $H= \frac 12trS,\
\tau=$det$S$. $H,\tau$  are  called  respectively an   affine   mean
and Gauss curvatures of   the surface $(M,f)$. For the Blaschke structure
on a surface $(M,f)$ the Affine Theorema Egregium holds:
$$K_h=H+J  \tag E$$
where $K_h$ denotes the Gauss curvature of the metric $h$.

\par
Two surfaces $(M,f),(N,g)$ we  call  equiaffinely  equivalent
iff there exist a  diffeomorphism $\phi :M \rightarrow  N$  and  an  affine
transformation $A \in $ ASL(3) such that $g\circ \phi  = A\circ f$.

The semi-isothermal coordinates on a semi Riemannian surface $(M,h)$ are
the local coordinates $(x_1,x_2)$ on
$M$ such that $h(\p_1,\p_1)=\e e^{-2u},\
h(\p_2,\p_2)=\eta e^{-2u}$, $h(\p_1,\p_2)=0$ where $\p_i=\frac{\p}{\p
x_i},\
u\in C^{\infty}(M)$ and
$\e\in\{-1,1\}$. It means that
$$h=e^{-2u}(\e dx_1\otimes dx_1+\eta dx_2\otimes dx_2).$$
In the semi-isothermal coordinates the Christoffel symbols of the
Levi-Civita connection of $h$ are given by:
$$
\align
\Gamma ^{1}_{11} &= -\p_{1}u,
\ \Gamma ^{2}_{11}=\e\eta \p_{2}u\\
\Gamma ^{1}_{12} &= -\p_{2}u,
\ \Gamma ^{2}_{12} = -\p_{1}u  \tag S-I  \\
\Gamma ^{1}_{22} &=\e\eta\p_{1}u,
\ \Gamma^{2}_{22}=-\p_{2}u
\endalign
$$

The Laplacian $\Delta$ on the space $(M,h)$ is the differential operator
defined in local coordinates by the formula:
$$\Delta \phi=\Theta^{-1}\left( \Sigma\p_j(\Theta h^{ij}\p_i\phi)\right)$$
where $\Theta=\sqrt{\mid det h\mid}.$
In the semi-isothermal coordinates we have:
$$\Delta \phi = e^{2u}(\e \p_1^2\phi+\eta\p_2^2\phi)=e^{2u}\Delta_0\phi$$
where  $\Delta_0=\e\p_1^2+\eta\p_2^2$.  If $K_h$ is the Gauss curvature of
 the
semi-Riemannian surface $(M,h)$ (i.e. $R(X,Y)Z=K_h (h(Y,Z)X-h(X,Z)Y)$ where
$R$ is the curvature tensor of $(M,h)$) then
$$\Delta u=K_h \tag $\Delta$ $$
which means that $\Delta_0 u =e^{-2u}K_h$.

\par
\bigskip
{\bf 2.  Projectively flat Blaschke   structures induced on surfaces.}
Let  us  recall
that an equiaffine connection $\nabla $  is  projectively  flat  iff  the
following conditions hold (for $\gamma (X,Y) = \frac 1{n-1}$ Ric$(X,Y) =
 \frac 1{n-1} tr\{Z\rightarrow R(Z,X)Y\})$
$$ \gather
 R(X,Y)Z = \gamma (Y,Z)X - \gamma (X,Z)Y \tag P1
\\ \nabla \gamma (X,Y,Z) = \nabla \gamma (Y,X,Z). \tag P2
\endgather $$
\noindent If $M$ is a surface then $(P1)$ always holds. If $\dim  M > 2$  then $(P2)$
is a consequence of $(P1)$ (see \cite{N-P-1}). We say that an  equiaffine
structure $(\nabla ,h,S)$ is projectively flat if $\nabla $ is a projectively flat
connection. The following lemma gives a description of  equiaffine
projectively flat structures induced by an affine immersion:
\par
\bigskip
{\bf Lemma 1.} {\it Let} $(M,f)$ {\it be  a  nondegenerate surface  in}
$\Bbb{R}^{3}$ {\it and} $(\nabla ,h,S)$ {\it be an equiaffine structure induced by} $f$  {\it and  an
equiaffine  transversal  field} $\xi ${\it . Then the induced connection }
 $\nabla $   {\it is projectively flat if and only if} $tr_{h}\nabla S = 0.$
\par
\medskip
{\bf Proof:}  As $M$  is
a  surface  we  have
$$R(X,Y)Z = \gamma (Y,Z)X - \gamma (X,Z)Y. \tag2.1$$
Hence we obtain
$$
 \nabla _{W}R(X,Y)Z = \nabla _{W}\gamma (Y,Z)X - \nabla _{W}\gamma (X,Z)Y.  \tag 2.2
$$
From (2.2) it follows that $\nabla $ is projectively flat iff
$$ \nabla _{W}R(X,Y)Z = \nabla _{Z}R(X,Y)W \tag 2.3 $$
for all $X,Y,W$,Z. On the other hand  from  the  Gauss  equation  we
obtain:
$$\gather
\nabla_{W}R(X,Y)Z = C(W,Y,Z)SX - C(W,X,Z)SY +  \tag 2.4
\\ h(Y,Z)\nabla S(W,X)- h(X,Z) \nabla S(W,Y).  \endgather $$
Hence from (Codazzi 1) and (2.3) it follows
$$\gather
 h(Y,Z)\nabla S(W,X) - h(X,Z)\nabla S(W,Y) \tag 2.5\\ = h(Y,W)\nabla S(Z,X) - h(X,W)\nabla S(Z,Y).
\endgather $$
Let $\{E_{1},E_{2}\}$ be an orthonormal frame, $h(E_{i},E_{i}) = \epsilon _{i} \in  \{-1,1\}$,  and
take $X=W=E_{1}, Z=Y=E_{2}$. We get
$$
 \epsilon _{1}\nabla S(E_{1},E_{1}) + \epsilon _{2}\nabla S(E_{2},E_{2}) = 0
\tag 2.6$$
which means $tr_{h}\nabla S = 0.$ It is not difficult to see that if the last
equation is satisfied then  (2.5)  holds  and  consequently $\nabla $  is
projectively flat.
$\diamondsuit $
\par
\medskip

Our present aim is to describe projectively
flat Blaschke connections $\nabla $  on  surfaces  as  well  as  Blaschke
immersions $f$ inducing such connections such that  affine Gauss  curvature
$\tau$ of $(M,f)$  is  constant. By $\Bbb{C}$-diagonalizable endomorphism
$S$ we mean an endomorphizm whose complexification is diagonalizable.
\par
\medskip
{\bf Theorem 1. } {\it Let} $f:M \rightarrow   \Bbb{R}^{3}$  {\it be  an  affine  immersion
with an equiaffine structure} $(\nabla ,h,S)$ {\it inducing a projectively  flat
Blaschke connection} $\nabla $ {\it and such that} $\tau = \det  S$ {\it is constant. If} $S$ {\it is}
$\Bbb{C}$-{\it diagonalizable then} $(M,f)$  {\it is  affinely  equivalent  to a locally
symmetric affine surface or to an affine sphere.}
\par
\medskip
{\bf Proof:} Let us recall that:
$$\gamma (X,Y) = 2H h(X,Y) - h(X,SY).\tag 2.7$$
\noindent Hence
$$\gather\nabla _{Z}\gamma (X,Y) = 2ZH h(X,Y) + 2H C(Z,X,Y)
- C(Z,X,SY) \tag 2.8\\
 - h(X,\nabla _{Z}S(Y)). \endgather$$
\noindent From (2.8) and Codazzi equations it follows  that $\nabla \gamma $  is  totally
symmetric iff
\par
$$
 - C(Z,X,SY) + 2ZH h(X,Y) = - C(Y,X,SZ) + 2YH h(X,Z)
\tag 2.9$$
\noindent or equivalently if
\par
$$
 2K(Z,SY) -2K(SZ,Y) = 2((YH)Z - (ZH)Y).
\tag 2.10$$
\noindent Hence for a projectively flat structure we get

$$\hn S(Z,Y) - \hn S(Y,Z) = (ZH)Y - (YH)Z. \tag 2.11$$
\par
\medskip
Next we consider two cases.
\par
\medskip
(a) The shape operator $S$ is $\Bbb{R}$-diagonalizable on $M$. Let $\lambda ,\mu $ be
eigenvalues of $S$ and $U = \{ x:\lambda (x) \neq  \mu (x)\}$.  Let  us  note  that
(int$(M\backslash U),f)$ is an affine sphere. We shall show that $(U,f)$  is  a
locally  symmetric  affine  surface.  Let $\{E_{1},E_{2}\}$  be   a   local
orthonormal frame, $h(E_{1},E_{1}) = \e, h(E_{2},E_{2}) = \eta  \in \{-1,1\}$ such that
\par
$$
 SE_{1} = \lambda E_{1}\qquad SE_{2} = \mu E_{2}.
\tag 2.12$$
We have $\hn_{X}E_{1} =\e \omega (X)E_{2}$,
$\hn_{X}E_{2} = - \eta \omega (X)E_{1}$  where $\omega ^{2}_{1} =\e\omega$  is  a
connection form for $\nabla $. From (2.12) we obtain
$$
\hn S(E_{2},E_{1}) = (E_{2}\lambda )E_{1} + \e(\lambda  - \mu )\omega (E_{2})E_{2}.
\tag 2.13$$
\noindent Analogously we obtain:
$$
\hn S(E_{1},E_{2}) = \eta (\lambda  - \mu )\omega (E_{1})E_{1}+ (E_{1}\mu )E_{2}.
\tag 2.14$$
\noindent Hence we get
$$ \gather
\hn S(E_{1},E_{2}) - \hn S(E_{2},E_{1}) =
(\eta (\lambda  - \mu )\omega(E_{1}) - (E_{2}\lambda ))E_{1} +
\tag 2.15$$
\\((E_{1}\mu ) - \e(\lambda  - \mu )\omega (E_{2}))E_{2}. \endgather$$
\noindent From (2.11) we have also
\par
$$
\hn S(E_{1},E_{2}) -\hn S(E_{2},E_{1}) = (E_{1}H)E_{2} - (E_{2}H)E_{1}.
\tag 2.16$$
\noindent Hence comparing (2.15) and (2.16) we obtain:
$$
- \eta (\mu  - \lambda )\omega (E_{1}) = \frac12 E_{2}(\lambda  - \mu ),
\ \ \e(\lambda  - \mu )\omega (E_{2}) = - \frac12 E_{1}(\lambda  - \mu )
$$
\noindent and consequently
$$
 \eta \omega (E_{1}) = \frac 12 E_{2}(\ln \mid \lambda  - \mu \mid ),
 \ \ \e\omega (E_{2}) = - \frac12 E_{1}(\ln \mid \lambda  - \mu \mid ).
\tag 2.17$$
Hence
$$\hn_{E_{1}}E_{2} = - \frac 12 E_{2}(\ln \mid \lambda  - \mu \mid )E_{1},\hn_{E_{2}}E_{1} = - \frac12 E_{1}(\ln \mid \lambda  - \mu \mid )E_{2}.
\tag 2.18$$
\noindent Let us introduce coordinates $(x_{1},x_{2})$ for which $E_{1} = \phi \partial _{1}, E_{2} = \psi \partial _{2}$
for some smooth functions $\phi ,\psi $. From  (2.18)  we  obtain  in  those
coordinates:
$$\gather
\phi \partial _{1}\psi \partial _{2} +
\phi \psi \hn_{\partial _{1}}\partial _{2} =
- \frac12 \phi \psi  \partial _{2}(\ln \mid \lambda  - \mu \mid )
\partial _{1}
\\ \psi \partial _{2}\phi \partial _{1} +
\phi \psi \hn_{\partial _{2}}\partial _{1} = - \frac 12 \phi \psi
\partial _{1}(\ln \mid \lambda  - \mu \mid )\partial _{2}.
\endgather$$
As $\hn$  is without torsion it yields:
$$
- \partial _{1}\ln \mid \psi \mid \partial _{2} - \frac 12
\partial _{2}(\ln \mid \lambda  - \mu \mid )\partial _{1} =
- \partial _{2}\ln \mid \phi \mid \partial _{1} -
\frac12 \partial _{1}(\ln \mid \lambda  - \mu \mid )\partial _{2}.
$$
\noindent Hence
$$ \partial _{2}(\ln \frac {\mid \phi \mid}{\sqrt{\mid \lambda -\mu \mid}}) = 0,
 \partial _{1}(\ln \frac{\mid \psi \mid}{\sqrt{\mid \lambda -\mu \mid }} )
  = 0. \tag 2.19 $$
\noindent From (2.19) it follows $\psi  = \beta (x_{2})\sqrt{\mid \lambda-\mu \mid}$ and $\phi  = \alpha (x_{1})\sqrt{\mid \lambda-\mu \mid}  $. Let us
introduce new coordinates $y_{1},y_{2}$ such that:
$$
y_{1} = \int\frac 1{\alpha}(x_{1}), y_{2} =\int \frac1{\beta}(x_{2}).
$$
\noindent It is clear that in new coordinates $E_{1} =
\sqrt{\mid \lambda- \mu\mid} \partial _{1}, E_{2} =
\sqrt{\mid\lambda-\mu\mid}\partial _{2}$
\par
\noindent and $h_{11} = \frac {\e}{\mid \lambda - \mu \mid } $ ,
$ h_{22} =\frac \eta{ \mid \lambda - \mu \mid } $.
From (2.10) we get
$$(\lambda  - \mu )K(E_{1},E_{2}) = (E_{1}H)E_{2} - (E_{2}H)E_{1}.\tag 2.20$$
By (2.20) it is clear that in introduced coordinates the
following equations are satisfied:
$$\align
K(\p_{1},\p_{2})&=\frac1{\lambda-\mu}(-\p_{2}H\p_{1}+\p_{1}H\p_{2})\\
K(\p_{1},\p_{1})&=\frac1{\lambda-\mu}(-\p_{1}H\p_{1}-\e\eta \p_{2}H\p_{2}) \tag K\\
K(\p_{2},\p_{2})&=\frac1{\lambda-\mu}(\e\eta\p_{1}H\p_{1}+\p_{2}H\p_{2}) \\
h_{11}&=\frac{\e}{\mid\lambda-\mu\mid} ,h_{22}=\frac{\eta}{\mid\lambda-\mu\mid} \tag h\\
S\p_{1}&=\lambda\p_{1},S\p_{2}=\mu\p_{2}  \tag S
\endalign
$$
By (K) and (h) one can easily obtain (using equality $\n=\hn+K$) the
expressions for connection coefficients of $\n$ in a chart $(x_{1},x_{2})$:
\par
\medskip
$$
\align
\Gamma ^{1}_{11} &= -\frac{\p_{1}\lambda}{\lambda-\mu},
\ \Gamma ^{2}_{11}=-\frac{\e\eta\p_{2}\mu}{\lambda-\mu}\\
\Gamma ^{1}_{12} &= -\frac{\p_{2}\lambda}{\lambda-\mu},
\ \Gamma ^{2}_{12} = \frac{\p_{1}\mu}{\lambda-\mu}  \tag $\Gamma$  \\
\Gamma ^{1}_{22} &=\frac{\e\eta\p_{1}\lambda}{\lambda-\mu},
\ \Gamma^{2}_{22}=\frac{\p_{2}\mu}{\lambda-\mu}
\endalign
$$
On the other hand if we define a structure $\str$ by
($\Gamma$), (h)  and (S)
then (K) holds, where $\n=\hn+K$ and $\hn$ is the Levi-Civita connection
for $h$,
and the Codazzi and Ricci equations are satisfied. For example
the Codazzi equation $\n S(X,Y)=\n S(Y,X)$ is equivalent to
$$\Gamma ^{1}_{12} = -\frac{\p_{2}\lambda}{\lambda-\mu},
  \Gamma ^{2}_{12} = \frac{\p_{1}\mu}{\lambda-\mu}$$
(see \cite{J-1}). It follows that $\str$ is an induced Blaschke structure of
a certain projectively flat affine surface with diagonalizable shape operator
in $\Bbb{R}^{3}$ iff the Gauss equation is satisfied.
It is easy to check that the Gauss equation is equivalent to the following
system of nonlinear partial differential equations of a second order
$$
\align
\Delta_{0}\mu + \frac2{\lambda-\mu}\mid\n\mu\mid^{2} &=\a\mu\\
\Delta_{0}\lambda - \frac2{\lambda-\mu}\mid \n \lambda \mid^{2} &
=-\a\lambda  \tag G \\
\p_{1}\lambda\p_{2}\mu= \p_{1}\mu\p_{2}\lambda
\endalign
$$
where $ \Delta_{0}=\e\p_{1}^{2}+\eta\p_{2}^{2}, \
\mid \n f\mid^{2}=\e \p_{1}f^{2}+ \eta\p_{2}f^{2},\
\a=\text{sgn}(\mu-\lb)
\in\{-1,1\}$. We can assume that $\a=1$ in another case changing
$\e$ by $-\e$ and $\eta$ by $-\eta$. So there is one-to-one
correspondence between  solutions of (G) and projectively flat surfaces
with diagonalizable shape operator. If $\lambda \mu  = \tau$ is constant,
then equation (G) can be reformulated after some simple computations as
(we set $\a=1$)
\par
$$
 \Delta _{0} \lambda  + \frac{2\lambda}{ \tau-\lambda ^{2}}
 (\e (\partial _{1}\lambda )^{2} + \eta (\partial _{2}\lambda )^{2})
 =-\lambda . \tag 2.21 $$
\noindent
 Next we consider three cases: $\tau>0,\ \tau=0,\ \tau<0$. Note
 that equation (2.21) is hyperbolic if and only if $\e\eta=-1$.
 \par
 \bigskip
 (i) $\tau>0$ .  As $\lb\ne\mu$ on $U$  and $\lb\mu=\tau$ we can assume
 that for example $\mid\lb\mid<\sqrt{\tau}$. Let us define
 a function $\Psi$ by
 $\Psi=-2\text{arc tgh}(\frac{\lb}{\sqrt{\tau}})$. Then
 equation (2.21) is:
 $$\Delta_0\Psi=-\text{sinh}\Psi.\tag$\Delta_1$ $$
 \par
 \bigskip
 (ii) $\tau=0$.  In that case if we define $\Psi=\frac1{\lb}$ where $\lb$ is
 the nonzero eigenvalue of $S$ then equation (2.21) is:
 $$\Delta_0\Psi=-\Psi.\tag$\Delta_2$ $$
 \par
 \bigskip
(iii) $\tau<0$ . Let us define a function $\Psi$ by
 $\Psi=-2\text{arc ctg}(\frac{\lb}{\sqrt{-\tau}})$. Then
 equation (2.21) is:
 $$\Delta_0\Psi=\text{sin}\Psi.\tag$\Delta_3$ $$

 In every case above the surface $(M,f)$ with constant $\tau$
 is an affine locally symmetric
 surface  which  follows
from [J-1] and the uniqueness theorem in affine differential geometry (see
\cite{D-N-V}). These surfaces have constant Gauss curvature with
respect to an appropriate nondegenerate scalar product in $\Bbb{R}^3$
(definite or indefinite).
\par
\bigskip
(b) Now let us assume that $(\nabla ,h,S)$  is  projectively  flat  and
that a shape operator $S$ has a  complex  eigenvalue.  In \cite{J-1}  to
describe locally symmetric surfaces with such a  shape  operator  we
used only assumption $tr_{h}\nabla S = 0$ and $\tau$ is constant which holds  also
for projectively flat structures with a constant curvature.  Hence
one can repeat the proof literally for  our  case.  In  particular
equation $(E4)$ is also a consequence  of  the  Gauss  equation  for
$(\nabla ,h,S)$ or of affine Theorema Egregium. Consequently in that  case
a projectively flat connection is locally symmetric. Let us note that for a
general projectively flat surface in the considered case the considerations
in  \cite{J-1} are still valid and we have in
coordinates introduced in \cite{J-1} using notation introduced there
(we take $b'=b, c'=-b$ see \cite{J-1}, p.217 as we can assume $\alpha=1$)
$$\align
S\p_{1} &= a\p_{1}+b\p_{2},\ S\p_{2} = -b\p_{1}+a\p_{2}
 \tag S'\\
h_{12} &= \frac1{b}, \ h_{ii}=0   \tag h'\\
\Gamma ^{1}_{11} &= -\frac{\p_{1}b}{b},
\ \Gamma ^{2}_{11}=\frac{\p_{1}a}{b}\\
\Gamma ^{1}_{12} &= 0,
\ \Gamma ^{2}_{12} =0 \tag $\Gamma$'  \\
\Gamma ^{1}_{22} &=-\frac{\p_{2}a}{b},
\ \Gamma^{2}_{22}=-\frac{\p_{2}b}{b}
\endalign
$$
and a structure $\str$ defined by these equations is integrable iff
$a,b$ satisfy the following equations (the Gauss equation)

$$
b\p_{12}a -2\p_{1}a\p_{2}b = -b^{2},
b\p_{12}b +\p_{1}a\p_{2}a-\p_{1}b\p_{2}b=ab
\   \tag G'$$
\par
\medskip
If the affine Gauss curvature $\tau$ of the above surface is constant then
$a^2+b^2=\tau=$ const and let us define a function $\phi$ by the equations:
$$a=\sqrt{\tau}\sin \phi,\ b=\sqrt{\tau}\cos \phi.$$
Then equation (G') is equivalent to
$$\text{sin}\phi\p^2_{12}\phi-\p_2\phi\p_2\phi\text{cos}\phi=\sin \phi.
 \tag 2.22$$
Let us define $\Psi=2 \text{arc tg}e^{\phi}$. Then equation (2.22) is
$$\p_1\p_2 \Psi=\text{cosh}\Psi.  \tag $\Delta_4$ $$
Every surface with constant $\tau$ considered  in that case is
a semi-Rieman\-nian surface with constant Gauss curvture
immersed in the Lorentz space $\Bbb{R}^3$.

It is known  that affine locally symmetric  surfaces  with
nondegenerate $\gamma $ have constant Gauss curvature with respect to some
scalar product in $\Bbb{R}^{3}$ and that this induced semi-Riemannian
structure coincides with the Blaschke structure.
It follows easily from the fact, that  affine
normal $\xi $ of such  surfaces  lies  on  a  centro-affine
  quadric $\Sigma
=\{ X: L(X,X) = 1 \}$,  where $L$ is a nondegenerate
symmetric form on $ \Bbb{R}^{3}$
(see \cite{J-1}). Differentiating an equation $L(\xi,\xi)=1$ we get for every
$X \in TM$ $ L(\xi_{*}(X),\xi) =0$ and consequently $L(f_{*}(SX),\xi)=0$. Hence
$\xi$ is a normal Riemanian field for $(M,f)$ with respect to $L$ and $(M,f)$ is
 nondegenerate  submanifold of $(\Bbb{R}^{3}, L)$. $L$ is a scalar product with
 respect to which $(M,f)$ has a constant nonzero Gauss curvature as
 the Riemannian structure coincides with the Blaschke structure, in particular is
 locally symmetric.
On the other  hand  every  surface  with  such  a
property  is  affine  locally  symmetric.  In fact if $g$ is the induced
semi-Riemannian metric then for the second fundamental form we have
$h(X,Y)=g(SX,Y)$ where $S$ is the (Riemannian) shape operator. Hence,
where by $\nu_h,\ \nu_g$ we denote volume forms of metrics $h$ and $g$,
$\nu_h=\sqrt{\mid \det S\mid}\nu_g$ and if $\n g=0$ then
$\n\nu_h=\sqrt{\mid\det S\mid }\n\nu_g=0$ as $\det S$ is constant.
The last equation is equivalent to the apolarity condition (A).
Hence the induced semi-Riemannian connection coincides with the
induced Blaschke connection.
\par
\bigskip

{\bf 3.  Affine locally strongly convex  spheres.}
A surface $(M,f)$ is called an {\it affine sphere} if $S  =  HI$  for
some $H \in \Bbb{R}$. A surface
$(M,f)$ is an affine sphere if and only if  $\nh K$ is a symmetric tensor.  The
other condition characterizing uniquely an affine  sphere  is  the
equation
$$
 R(X,Y)Z = \hat R(X,Y)Z + [ K_{X},K_{Y} ]Z \tag R2
$$
  Let $(M,f)$  be  an
affine locally strongly convex surface in $\Bbb{R}^{3}$  with  a  Blaschke
structure $\str$. Then the following lemmas hold.
\par
\bigskip
{\bf Lemma A.} {\it Let} $x_{0} \in  M$ {\it satisfies the condition}
$K_{x_{0}}\neq  0.$ {\it Then  there
exists an open neighborhood} $V$ {\it of} $x_{0}$
{\it and a local orthonormal  frame}
$\{ X,Y \}$ {\it defined on} $V$ {\it and satisfying the equations:}
$$
K(X,X) = - \lambda Y, K(X,Y) = - \lambda X, K(Y,Y) = \lambda Y. \tag 3.1
$$
\par
\noindent {\it with} $\lambda  = \frac12 \mid K \mid= \frac 14 \sqrt{h(C,C)}$.
{\it In the basis} $\{X,Y\}$ {\it the endomorphisms}
$K_{X},K_{Y}$ {\it are represented by the following matrices} :
\par
$$
 K_{X} = \pmatrix 0&-\lambda\\ \cr-\lambda &0\endpmatrix \
 K_{Y} = \pmatrix -\lambda &0 \\ \cr0&\lambda \endpmatrix  \tag 3.2
$$
\noindent {\it An acute angle between any two null directions of the  cubic
  form is} $\frac13 \pi $  {\it and  for  any  vectors} $U,V,W$
   {\it the  endomorphisms} $K_{U},K_{V}$
{\it satisfy the equation}
\par
$$
 [ K_{U},K_{V} ] (W) = - J(h(V,W)U - h(U,W)V) \tag 3.3
$$
{\bf Proof:}  Let $\{ E_{1},E_{2}\}$ be  an  orthonormal
local frame defined on an open set $V \subseteq  M$ and $x_{0} \in $ V.
 From ($R1)$  it
follows that there exist smooth functions
$a,b \in  C^{\infty }(V)$ such that
$$
 K(E_{1},E_{1}) = aE_{1} + bE_{2}, K(E_{1},E_{2}) = bE_{1} - aE_{2} \tag 3.4
$$
\noindent Let us take $X = \sin (\phi )E_{1} - \cos (\phi )E_{2},Y = \cos (\phi ) E_{1} + \sin (\phi ) E_{2}$
where $\phi $ is a smooth function. It is clear that $a^{2} + b^{2} =
\frac 1{16} h(C,C)$.
Let us take $\lambda  = \frac 14 \sqrt {h(C,C)}$. Then $J=2\lb^2$.
Locally  there  exists  a
function $\psi  \in  C^{\infty }(V)$ such that $a = \lambda \cos  \psi ,
b = \lambda \sin  \psi $. Notice that:
$$
\gather
K(Y,Y) = \\
(\cos (\phi )^{2} - \sin (\phi )^{2})K(E_{1},E_{1}) +
 2\sin (\phi )\cos (\phi )K(E_{1},E_{2}) \\
 = \cos (2\phi )K(E_{1},E_{1}) + \sin (2\phi )K(E_{1},E_{2})\\ =
 (a\cos (2\phi )
  +  b\sin (2\phi ))E_{1} +(b\cos (2\phi ) - a \sin (2\phi ))E_{2} \\=
\lambda ( \cos (2\phi  - \psi )E_{1} + \sin (\psi  - 2\phi )E_{2}) \\=
\lambda (\cos (\psi  - 2\phi )E_{1} + \sin (\psi  - 2\phi )E_{2}).
\endgather
$$
\par
Hence $K(Y,Y) = \lambda Y$ if and only if
$\phi  = \frac 13 \psi  + \frac {2k\pi} 3,k \in  \Bbb{Z}$. It is easy
to check using $(R1)$ that  with  such  a  choice  of $\phi $  the  other
equations are satisfied. Let us note  that $X$  lies  on  the  null
direction of $C$ (see \cite{N-P-3}). As the angles $\phi $ and
$\phi  + \pi $  give  the
same null direction  of $C$  the  lemma is proved.$\diamondsuit$
\par
\bigskip
{\bf Lemma B. }  {\it Let} $(M,f)$   {\it be   an
affine   locally   strongly convex surface with induced  Blaschke  structure  and}
 $\{ E_{1},E_{2} \}$  {\it be  a local  orthonormal  frame  satisfying}  (1).
 {\it   Then   the   following equations are satisfied :}
$$
\align
\nh K(X,E_{1},E_{1}) &= - 3\lambda  \omega (X)E_{1} - (X\lambda )E_{2} \\
\nh K(X,E_{1},E_{2}) &= - (X\lambda )E_{1} + 3\lambda  \omega (X)E_{2}
\tag 3.5 \\
\nh K(X,E_{2},E_{2}) &= 3\lambda  \omega (X)E_{1} + (X\lambda )E_{2}
\endalign
$$
\noindent {\it with} $\omega  = \omega ^{1}_{2}$
 {\it the connection form defined by  }
 $\nh_{X}E_{i} = \omega ^{j}_{i}(X)E_{j}$.
\par
\medskip
{\bf Corollary.} {\it Let} $(M,f)$ {\it be an affine
locally strongly  convex sphere in} $ \Bbb{R}^{3}$.
{\it Let us define} $U := \{ x \in  M: h(C,C) > 0 \}$.  {\it Then  for
every} $x_{0} \in  U$ {\it there exists a  local  coordinate  system}
$(V,x_{1},x_{2})${\it such that} $V$ {\it is a neighborhood of}
$x_{0}, V \subseteq  U${\it , and an equation}
\par
$$
 E_{1} = e^{u} \partial _{1} , E_{2} = e^{u} \partial _{2} \tag{E1}
$$
\noindent {\it holds, where}
$u = \frac13\ln \lambda ,\lambda  = \frac 14 \sqrt{h(C,C)}$ {\it and}
 $\{ E_{1},E_{2} \}$  {\it is  a  local frame satisfying  equations} (1).
  {\it The coordinates} $(V,x_{1},x_{2})$  {\it are
isothermal coordinates for} $(M,h)$ {\it and the following  equations  are
satisfied:}
\par
$$
h(\partial _{1},\partial _{1}) = h(\partial _{2},\partial _{2}) = e^{-2u},
h=e^{-2u}((dx_{1})^{2}+(dx_{2})^{2})
\tag h
$$
$$\gather
 K(\partial _{1},\partial _{1}) = - e^{2u} \partial _{2} ,
 \  K(\partial _{2},\partial _{2}) = e^{2u} \partial _{2},\tag K
\\ K(\partial _{1},\partial _{2})= - e^{2u} \partial _{1}
\endgather
$$
\noindent {\it and}
$$
\align
\nabla _{\partial _{1}}\partial _{1} &= - \partial _{1}u \partial _{1} +
 ( \partial _{2}u - e^{2u}) \partial _{2}\\
 \nabla _{\partial _{1}}\partial _{2} &= - ( \partial _{2}u + e^{2u})
  \partial _{1} - \partial _{1}u \partial _{2} \tag $\nabla$ \\
\nabla _{\partial _{2}}\partial _{2} &= \partial _{1}u \partial _{1}
 - ( \partial _{2}u - e^{2u}) \partial _{2}
\endalign
$$
$$
 \Delta _{0} u = e^{-2u}(H + 2e^{6u}). \tag $\Lambda$
$$
\noindent {\it or equivalently}
\par
$$
 \Delta _{0} \ln  \lambda  = 3 \lambda  ^{-\frac 23}
  (H + 2\lambda ^{ 2}) \tag $\Lambda'$
$$
\noindent {\it where} $\Delta _{0}$ {\it is the standard Laplacian on}
$\Bbb{R}^{2}$. {\it The condition} $h(C,C)$ = const
({\it or equivalently} $K_{h}$ = const. ) {\it  is equivalent to }
$\nh K = 0 $
(  {\it or equivalently} $\nh C = 0$ ).
\par
\medskip
{\bf Proof:}  From (3.5) we obtain that $\nh K$ is symmetric if  and  only
if  equations $(\omega )$ below are satisfied.
\par
$$
 \omega (E_{1}) = - du(E_{2}) ,\ \omega (E_{2}) = du(E_{1})  \tag $\omega$
$$
Let $(W,y_{1},y_{2})$ be a coordinate system such that
\par
$$
E_{1} = \phi  \partial _{1} ,\  E_{2} = \psi  \partial _{2}
$$
\noindent for some smooth functions $\phi ,\psi $. We can assume that $\phi ,\psi $ are positive
in an opposite case  changing   coordinates  as  follows  :$(y_{1},y_{2})
\rightarrow  (\epsilon y_{1},\eta y_{2})$ where $\epsilon ,\eta  \in  \{ -1,1 \}$. The equations
$$ \nh_{E_{2}} E_{1} = - \omega (E_{2}) E_{2},\
\nh_{E_{1}} E_{2} = \omega (E_{1}) E_{1}
\tag 3.6$$
\noindent hold on U. Hence we get
\par
$$
\gather
\partial _{1}\psi  \partial _{2} +
\psi \nh_{\partial _{1}}\partial _{2} =
- \psi  \partial _{2}u \partial _{1}\\
\partial _{2}\phi  \partial _{1} +
\phi \nh_{\partial _{2}}\partial _{1} =
- \phi  \partial _{1}u \partial _{2}.
\endgather
$$
\noindent As $\nh$  is without torsion we obtain
$$
 \partial _{1}\ln  \psi  = \partial _{1}u\qquad \partial _{2}\ln  \phi
 = \partial _{2}u.  \tag 3.7
$$
\noindent From (3.7) it follows that there exist smooth positive functions
$\Phi ,\Psi $ such that
\par
$$
 \phi  = e^{u}\Phi (y_{1})\qquad \psi  = e^{u}\Psi (y_{2}) \tag 3.8
$$
\noindent Let us change coordinates in the following way
\par
$$
x_{1} = \int (\frac1{\Phi}) (y_{1}) \qquad x_{2} = \int (\frac1{\Psi}) (y_{2})
$$
\noindent It is clear that in new coordinates equations $(E1)$  hold.  Let
us note also that from $3), (R2)$ and  equations
\par
$$
\gather
R(X,Y)Z = H(h(Y,Z)X - h(X,Z)Y) \tag 3.9 \\
\hat R(X,Y)Z = K_{h}(h(Y,Z)X - h(X,Z)Y) \tag 3.10
\endgather
$$
\noindent where $H$ is an affine mean curvature of an affine sphere
$(M,f)$  and
$K_{h}$ is the Gauss curvature of $(M,h)$, we get the
 Affine Theorema Egregium
$$
 K_{h} = H + 2 e^{6u} = H + \frac18 h(C,C) \tag E
$$
Equation (E) is valid for any  affine  surface  as
one can check computing the Ricci tensor  of  the  affine  metric.
From (E) it is clear that $K_{h}$ = const iff $h(C,C)$ = const. It  also
follows from (3.5) and $(\omega )$  that $h(C,C)$ is constant iff $\nh K = 0.$ In  a
case $h(C,C) = 0$ it is obvious as then $K = 0.$ Let us note also that
in new coordinates
\par
$$
 h = e^{-2u}(( dx_{1})^{2} + (dx_{2})^{2}) \tag 3.11
$$
\noindent Hence
$$
 \Delta  u = K_{h} = H + 2e^{6u} \tag 3.12
$$
\noindent where $\Delta $ is the Laplacian for $(M,h)$. Equation (3.12) is
equivalent to
\par
$$
\Delta _{0} u = e^{-2u}(H + 2e^{6u}), H \in \Bbb{R}  \tag $\Lambda$
$$
\noindent where $\Delta _{0} = \partial ^{2}_{1} + \partial ^{2}_{2}$
is the standard Laplacian on $\Bbb{R}^{2}.\diamondsuit$
\par
\medskip
{\bf Theorem  2.}  {\it Let} $u$ {\it satisfies  equation} $(\Lambda )$ {\it on an  open  simply
connected set} $\Omega  \subseteq \Bbb{R}^{2}${\it .
Then there exists an immersion} $f_{u}: \Omega  \rightarrow \Bbb{R}^{3}$
{\it such that} $(\Omega ,f_{u})$ {\it is a locally strongly convex affine
  sphere  with an
affine mean curvature} $H$ {\it and Fubini-Pick invariant} $16e^{6u}$
{\it .  In  the
standard coordinates on} $\Omega $ {\it  equations}
(h),(K),$(\nabla )$ {\it hold where} $h$
{\it is an affine metric,} $K$ {\it is the difference tensor and}
$\nabla = \nh  + K$  {\it is the Blaschke connection for} $
(\Omega ,f_{u})${\it . An immersion} $f$ {\it is unique up  to
an affine transformation.}
\par
{\it Let} $(M,f)$  {\it be  a  locally  strongly  convex  affine
surface with induced Blaschke  structure which is an affine sphere  with
 affine mean curvature} $H${\it . Then} $M = U \cup \bar{V} = \bar{U} \cup  V$
 {\it where} $U,V$  {\it are
open subsets of} $M$ {\it such that} $U \cap  V  =   \emptyset$   {\it ,} $(U,f)$  {\it is  a  locally
strongly convex quadric and} $V = \{ x:h(C,C) > 0 \}$ . {\it The  function}
$\lambda:= \frac 14 \sqrt{h(C,C)}$ {\it satisfies on} $V$ {\it the equation}
\par
$$
\Delta  \ln  \lambda  = 3( H + 2\lambda ^{2})
$$
\noindent {\it and  around  any  point} $x_{0} \in  V$  {\it there
exist  local  isothermal
coordinates} $(W,x_{1},x_{2})$  {\it such  that  the  function}
 $u = \frac13 \ln \lambda \circ \gamma ^{-1}$
{\it satisfies on} $\gamma (W)$ {\it  equation} $(\Lambda )$ {\it where}
 $\gamma (p) = (x_{1}(p),x_{2}(p))$  {\it and}
$(W,f)$ {\it is equiaffinely equivalent to  an affine sphere in}
$ R^3$  {\it given by} $(\gamma (W),f_{u})$.
\par
\medskip
{\bf Proof :}  Let a function $u$ satisfies on $\Omega $  equation $
(\Lambda )$. Let
us define an affine metric $h, $a  connection $\nabla $  and  a  difference
tensor $K$ by  formulas $(h),(\nabla )$ and (K). Let  us  define  also  a
connection $\nh$  as $\nabla  = \nh  + K$. It is clear that $\nh$  is  the
 Levi-Civita connection for $h$. If we define an orthonormal frame
 $\{E_{1},E_{2}\}$ by  equations $(E1)$ and by $(\omega )$ the differential
 form $\omega $ then it easy  to check  that  equations  (3.1)  hold
 for $\{E_{1},E_{2}\}$  and $\omega $  is  the connection form
 $\omega ^{1}_{2}$ with respect to $\{E_{1},E_{2}\}$. From (3.5) it follows
 that $\nh K$ is  a  symmetric  tensor. Hence    equation (R2)  is
satisfied where $R$ is the curvature tensor of $\nabla $. Let us define a
shape operator by $S = H$ Id$_{T\Omega }$. From (R2), (3.10)  and  (3.3)  it  follows
that the Gauss equation (3.9) is equivalent to (E) which in  turn  is
equivalent to $(\Lambda )$. From (K) and an equation $\nabla = \nh+ K$  it  easily
follows that $C = \nabla h$ is a symmetric tensor. Hence  equations  of
Gauss,Codazzi and Ricci are satisfied for the  structure $(\nabla ,h,S)$.
From Radon's Theorem (see \cite{D-N-V}) it follows that there exists  a
nondegenerate immersion $f :\Omega  \rightarrow \Bbb{R}^{3}$ with $(\nabla ,h,S)$
 as the induced Blaschke  structure.  The   uniqueness   up   to   an   equiaffine
transformation follows from \cite{D}.
\par
The second part  of  the  theorem  follows  from  lemmas  and
Corollary.$\diamondsuit$
\par
\medskip
{\it Remark.}  Let us note that the equation ($\Lambda$)
is equivalent to the equation
$$\Delta_0\Psi=e^{2\Psi}+\e e^{-\Psi} \tag $\Delta_5$ $$
where $\e=\text{sgn}H\in\{-1,0,1\}$ and
$\Psi(x_1,x_2)=2u(\frac{x_1}a,\frac{x_2}a)-\ln b$ where
$\ \ a = (4\mid H\mid)^{\frac13}$, $b = (\frac{\mid H\mid}2)^{\frac13}$
if $H\ne 0$ and $\Psi=2u(\frac{x_1}2,\frac{x_2}2)$ if $H=0$.
\par
\medskip
It is interesting to know to what extent  an  affine  sphere  is
determined by its affine metric and the Fubini-Pick invariant. The
following corollary gives an answer to this question.
\par
\medskip
{\bf Corollary.}  {\it Let} $(M,f)$  {\it be  an  affine  sphere  in}
$\Bbb{R}^{3}$
  {\it with  a definite affine metric} $h$
  {\it and  Fubini-Pick invariant} $h(C,C)${\it . Let
us assume that} $h(C,C)$ {\it is positive on} M. {\it Then  for  every} $x_{0} \in  M$
{\it there exists a neighborhood} $V$ {\it of} $x_{0}$ {\it and a one parameter family  of
affine immersions} $\{f_{a}:a \in  O(2)\}$ {\it defined on} $V${\it , such that} $(V,f_{a})$  {\it is
an affine sphere and each immersion} $f_{a}$ {\it has the same induced affine
metric} $h$ {\it and the same Fubini-Pick invariant} $h(C,C)${\it .  Every  affine
sphere immersion whose  induced Blaschke structure has  an  affine
metric} $h$  {\it and  the  Fubini-Pick  invariant} $h(C,C)$   {\it is   locally
equiaffinely equivalent to one of immersions} $f_{a}${\it .}
\par
\medskip
{\bf Proof:}  Let $(M,f),(M,f')$ be two affine spheres with the same
affine metrics and Fubini-Pick invariants. From  the  theorem  for
every point $x_{0} \in  M$ there exist local charts $(U,x_{1},x_{2}),(U',x_{1}',x_{2}')$
around $x_{0}$  such  that    equations  (K),(11)  and $(\nabla )$   hold
respectively for $u = \frac13\ln \lambda ,
u' = \frac13\ln \lambda '$ where
$\lambda  = \frac14\sqrt{h(C,C)}\circ \gamma ^{-1}, \lambda '=
\frac14\sqrt{h(C,C)}\circ \gamma^{,-1}$ . Let us define
 $\phi  = \gamma '\circ \gamma ^{-1} = (x'(x_{1},x_{2}),x'(x_{1},x_{2}))$.
Then $\lambda  = \lambda '\circ  \phi $ and $h_{ij} = h'_{ij}\circ \phi $.
 If we denote by $A^{i}_{j} = \frac{\partial x'_{i}}{\partial x_{j}}$
  then  the
transformation rules for $h$ gives us the following equations:
 $(A^{1}_{1})^{2} + (A^{2}_{1})^{2} = 1, (A^{2}_{2})^{2} + (A^{1}_{2})^{2} = 1, A^{1}_{1}A^{1}_{2} + A^{2}_{1}A^{2}_{2} =
 0.$ Hence $A^{1}_{1} = \cos  \alpha ,
A^{2}_{1} = \sin  \alpha , A^{2}_{2} = \epsilon  \cos  \alpha ,
A^{1}_{2} = - \epsilon  \sin  \alpha $ for $\epsilon  \in  \{ -1,1 \}$
and $\alpha $ a
smooth function. As $\partial _{i}A^{k}_{j} = \partial _{j}A^{k}_{i}$
it is easy  to  check  that $\alpha $  is
constant. We conclude that $\phi  = g$ where
\par
$$
 g = \pmatrix \cos\alpha &-\epsilon \sin \alpha &a \\ \sin \alpha
 &\epsilon \cos \alpha &b \\ 0&0&1 \endpmatrix \in  AO(2), \quad
 \alpha ,a,b \in  R
\tag g$$
\noindent is an affine orthogonal transformation of $ \Bbb{R}^{2}$. The  group $AO(2) =
E(2)$ is a symmetry group of  equation $(\Lambda )$. Let us note that the
difference tensor $K'$ of $(M,f')$ in the chart
$(V,x_{1},x_{2}),V$ = dom $\gamma$ $\cap $
dom $\gamma '$, is represented as follows: ( see Lemma A )
\par
$$
\align
K'(\partial _{1},\partial _{1})& = - e^{2u}( \sin  3\alpha  \partial _{1}
+ \epsilon  \cos  3\alpha  \partial _{2})\\
K'(\partial _{2},\partial _{2})& = e^{2u}( \sin  3\alpha  \partial _{1} +
\epsilon  \cos  3\alpha  \partial _{2})  \tag $\alpha$ \\
K'(\partial _{1},\partial _{2})& =
- e^{2u}(\epsilon  \cos  3\alpha  \partial _{1}
- \sin  3 \alpha  \partial _{2})
\endalign
$$
\par
\noindent and $h ,\nh$  are the same for $f$ and $f'$. $h$ is given by (3.11)
and  $\nh$   is
determined by h. On the other hand let us assume that $u$  satisfies
 equation $(\Lambda )$ on a simply connected subset $\Omega$ of
 $ \Bbb{R}^{2}$, $ \Omega  = \gamma (U)$ and $h$
is given by (3.11). Let us denote by $\nh$   the  Levi-Civita  connection
for h. It is easy to show that tensor $\nh K_{\alpha }$ where
$K_{\alpha }$ is  defined  by
 formula $(\alpha )$ is symmetric. It is also clear that
$K_{\alpha }$  satisfies (R1). Hence one can show as in the proof of
the Theorem that there
exists an immersion $f_{a}$  with  an  induced  affine  metric $h$  and
a difference tensor $K_{\alpha }$ , where $\nabla  = \nh + K_{\alpha }$ and
$$
 a = \pmatrix \cos 3\alpha &-\epsilon \sin 3\alpha \\
 \sin 3\alpha &\epsilon \cos 3\alpha \endpmatrix
\tag a
$$
\noindent is an element of the group $O(2)$. Let us note that
$K_{\alpha } = K_{\alpha +\frac{2\pi}3}$. It
is obvious from the  construction  of  coordinates  as  the  angle
between two null directions is $\frac{2\pi}3.\diamondsuit$
\par
\medskip
{\it Remark.}  Let us take $M = \Omega $ and $f = f_{u}$ where $u$ satisfies
equation $(\Lambda )$. We can construct as above a family of immersions
 $f_{a}$.
Let us note that $f_{a}$ are  not  equiaffinely  equivalent  as  affine
immersions $f_{a}: \Omega  \rightarrow \Bbb{R}^{3}$.  It  is  interesting  to  know  the
conditions under which $( \Omega , f_{a})$ is equivalent to
$(\Omega , f_{id})$  as  an
affine hypersurface. We have the following :
\par
\medskip
{\bf Proposition.} {\it Let a function} $u$ {\it satisfies on} $\Omega $
{\it  equation} $(\Lambda )$
{\it and let} $(\Omega ,f_{u})$ {\it be an affine sphere with  affine metric  and
Blaschke  connection given by formulas} (3.11) {\it and} $(\nabla )$. {\it Let} $(\Omega ,f_{a})$
 {\it be as above,} $f_{id} = f_{u}${\it .
Then} $(\Omega ,f_{a})$ {\it is affinely equivalent to} $(\Omega ,f_{u})$
 {\it iff} $u = u\circ g$, {\it where a is given by the formula} (a)
 {\it and} $g$ {\it is given by} (g)
{\it for some} $a,b \in R${\it .}
\par
\medskip
{\bf Proof:}   Let us assume that there exists a  diffeomorphism $\phi
:\Omega  \rightarrow  \Omega $ such that $f_{u}\circ \phi  = A\circ f_{a}$ for an equiaffine  transformation
$A$. Then $\phi $ is an isometry and an affine diffeomorphism with respect
to $\nabla $ and $\nabla _{\alpha } = \nh + K_{\alpha }$ , which means
$\nabla _{\phi _{*}(X)} \phi _{*}(Y) = \phi _{*}(\nabla _{\alpha X}Y)$  for
every  vector  fields $X,Y$  on $\Omega $.   Hence
$$\phi _{*}( K_{\alpha }(X,Y) ) =
K(\phi _{*}(X),\phi _{*}(Y)).$$
 It follows
$h_{\alpha }(C,C) \circ  \phi ^{-1} = h_{id}(C,C)$. Hence $\lambda  \circ  \phi ^{-1}
= \lambda $ and $u\circ \phi ^{-1}$ = u. As $\phi $ is an isometry
and $u\circ  \phi ^{-1} = u$ it  is  easy
to show as above that $\phi $ has a form of  (g). On the other hand if
 $u\circ g = u$ then one can check  exactly as in  the  Corollary  above  that
$(\Omega ,f_{a})$ is  affinely  equivalent  to $(\Omega ,f_{u})$  which  concludes  the
proof.$\diamondsuit$
\par
\medskip
{\it Remark.}  Let us note  that  if $u$  is  constant  then  all
hypersurfaces $( \Omega ,f_{a})$ are equiaffinely equivalent. In particular a
locally strongly convex affine sphere  with  constant  Fubini-Pick
invariant $h(C,C) \neq  0$ is characterized uniquely by $h(C,C)$.
 In fact
every affine sphere with constant nonzero Fubini-Pick invariant
and definite affine metric h is
equiaffinely equivalent to  an open part of the surface $(X^{2} - Y^{2})Z
= \pm  \frac1c$ where $c \neq  0$ depends only on $h(C,C)$. The  fundamental  system
of equations for $f$ is very simple for constant $\lambda  ($ see $(\nabla ) )$  and
it is easy to see that if we take a chart for which $E_{i} = \partial _{i}$ then

$$f(x_{1},x_{2}) = ( e^{-\lambda x_{2}} \cosh (\sqrt{3}\lambda x_{1})
,e^{-\lambda x_{2}} \sinh (\sqrt{3}\lambda x_{1}),
\pm  \frac1c e^{2\lambda x_{2}} )$$
where $c = \frac{3\sqrt{3}}{128} h(C,C)^{2}$ (see also \cite{M-N}) ). From our Proposition it
follows that the surface $M = \{ (X,Y,Z):(X^{2} - Y^{2})Z =
\pm  \frac1c \}$ is
equiaffinely homogeneous. It is easy to check that $( \Bbb{R}^{2},f)$  is  an
orbit of the point $(1,0, \pm  \frac1c )$  by  the  group $G =
\Bbb{R} \oplus \Bbb{R}$
  of equiaffine transformations
\par
$$
 \pmatrix
 e^{-a} \cosh b & e^{-a} \sinh b & 0 \\
 e^{-a} \sinh b & e^{-a}\cosh b & 0 \\
 0 & 0 & e^{2a}
 \endpmatrix  a,b \in \Bbb{R} \oplus \Bbb{R}
 $$
\par
\bigskip
{\bf 4.  Affine spheres with an indefinite affine  metric.}
Let $(M,f)$
be an affine surface with  an  induced  Blaschke  structure  whose
affine metric is indefinite. Let $x_{0} \in  M$ be  any  point  of $M$  and
$\{E_{1},E_{2}\}$ be a local orthonormal frame defined on a  neighborhood $U$
of the point $x_{0}$, such that $h(E_{1},E_{1}) = 1,h(E_{2},E_{2}) = - 1.$
\par
Let $K$ be a difference tensor of the induced structure.Then the
equations
$$
 K(E_{1},E_{1}) = aE_{1} + bE_{2},K(E_{1},E_{2}) = -(bE_{1} + aE_{2})
 \tag 4.1
$$
\noindent hold for some smooth functions $a,b \in  C^{\infty }(U)$.
The  above  equations
are   a   consequence   of $(R1)$.  Let  us note that
$h(K(E_{1},E_{1}),K(E_{1},E_{1})) = a^{2} - b^{2}$ and $h(C,C) = 16 ( a^{2} - b^{2})$. Let  us
denote $\epsilon $ := sign $(h(C,C))$ and $\lambda : = \sqrt{\epsilon (a^{2} - b^{2})}$ . Let us assume that
$\lambda (x_{0}) \neq  0$ and define $V := \{ x:\epsilon (x) = \epsilon (x_{0})\} \cap $  U. $V$  is  an  open
neighborhood of $x_{0}$. Now we consider three cases :
\par
\bigskip
i) $\epsilon (x_{0}) = 1$ which means $(a^{2} - b^{2}) > 0.$ Then there  exists
a function $\psi  \in  C^{\infty }(V)$ such that $a =
\lambda \cosh (\psi ), b = \lambda \sinh  (\psi )$.  Let us
define a new orthonormal frame $\{ X,Y \}$ by the equations
$$
 X = \cosh (\phi )E_{1} + \sinh (\phi )E_{2},
 Y = \sinh (\phi )E_{1} + \cosh (\phi )E_{2} \tag 4.2
$$
\noindent where $\phi $ is a smooth function. Then

$$
\gather
K(X,X) = (\cosh (\phi )^{2} + \sinh (\phi )^{2})K(E_{1},E_{1})\\
+ 2\sinh (\phi )\cosh (\phi )K(E_{1},E_{2}) =
 \cosh (2\phi )K(E_{1},E_{1}) \\ + \sinh (2\phi )K(E_{1},E_{2}) =
(a \cosh (2\phi ) -  b\sinh (2\phi ))E_{1}
 \\ + (b\cosh (2\phi ) - a\sinh (2\phi ))E_{2} =
\lambda ( \cosh (2\phi  - \psi )E_{1} \\+ \sinh (\psi  - 2\phi )E_{2})=
 \lambda (\cosh (\psi  - 2\phi )E_{1}\\ + \sinh (\psi  - 2\phi )E_{2}).
\endgather
$$
\par
Hence $K(X,X) = \lambda X$ if and only if $\phi  = \frac13 \psi $.
From (4.2) it is  clear
that $K(X,Y) = - \lambda Y$ if we choose $\phi $ as above.
\par
\vskip 0,5cm
ii)  $ \epsilon (x_{0}) = - 1.$  Then there exists a  function $\psi  \in  C^{\infty }(V)$
such that $a = \lambda \sinh (\psi )$, $b = \lambda \cosh  (\psi )$.  Let  us  define  a  new
orthonormal frame $\{ X,Y \}$ by  equations (4.2). It is easy to show
as above that
\par
$$
K(X,X) = \lambda (\sinh (\psi  - 2\phi )E_{1} +
\cosh (\psi  - 2\phi )E_{2})
$$
Hence in that case $K(X,X) = \lambda Y, K(X,Y) = - \lambda X $
for $\phi  = \frac13 \psi .$
\par
\vskip 0,5cm
iii$)$  Now we consider the case $\epsilon (x_{0}) = 0.$ Then, in general, we  can
not expect that $\epsilon (x) = 0$ in a certain neighborhood of $x_{0}$.  Let  us
assume that there exists a neighborhood $V$ of $x_{0}$ such that $\epsilon (x) = 0$
for $x \in $ V. Then in $V$ an equalitty $ a^{2} = b^{2}$ holds and:
\par
$$
K(E_{1},E_{1}) = aE_{1} \pm  aE_{2} ,
K(E_{1},E_{2}) = - (\pm  aE_{1} + aE_{2})
$$
\noindent and we can assume $a = b$ in an opposite case changing $E_{1}$  by $-E_{1}$.
Let us define local fields $X,Y$ by $X = (E_{1} + E_{2}), Y = \frac12(E_{1} - E_{2})$.
Then it is easy to check that
\par
$$
K(X,X) = 0, K(X,Y) = 0, K(Y,Y) = aX
$$
\noindent and
\par
$$
h(X,X) = h(Y,Y) = 0 , h(X,Y) = 1
$$
\noindent Hence we have proved:
\par
\bigskip
{\bf Lemma A1.} {\it Let} $(M,f)$ {\it be as above and} $x_{0} \in $ M. {\it Let us  define} $\alpha
:= \epsilon (x_{0})$ {\it and take} $M_{\alpha }$ := int$(\{ x:\epsilon (x) = \alpha  \} \cap  M )$. {\it If} $x_{0} \in  M_{\alpha }$  {\it then
there exists a local  frame} $\{ E_{1},E_{2}\}$  {\it defined  on  a
neighborhood} $V_{\alpha } \subseteq  M_{\alpha }$ {\it of} $x_{0}$
 {\it and satisfying  respectively  for} $\alpha  =
1,-1,0$ {\it the equations :}
$$
\gather
\alpha = 1 \qquad
 K(E_{1},E_{1}) = \lambda E_{1} ,
  K(E_{1},E_{2}) = - \lambda E_{2}, K(E_{2},E_{2}) = \lambda E_{1}.
  \tag 4.3 \\
h(E_{1},E_{1}) = 1, h(E_{1},E_{2}) = 0, h(E_{2},E_{2}) = - 1; \\
 \alpha  = - 1 \qquad K(E_{1},E_{1}) = \lambda E_{1} ,
 K(E_{1},E_{2}) = - \lambda E_{2}, K(E_{2},E_{2}) = \lambda E_{1}.
 \tag 4.4 \\
h(E_{1},E_{1}) = - 1, h(E_{1},E_{2}) = 0, h(E_{2},E_{2}) = 1; \\
\alpha  = 0 \qquad K(E_{1},E_{1}) = 0,\ \  K(E_{1},E_{2}) = 0,
\ \ \ K(E_{2},E_{2}) = aE_{1} \tag 4.5 \\
h(E_{1},E_{1}) = h(E_{2},E_{2}) = 0 , h(E_{1},E_{2}) = 1
\endgather
$$
\par
\medskip
\noindent {\it where} $\lambda  = \frac14\sqrt{\alpha h(C,C)}$
{\it and} $a \in  C^{\infty }(V_{0})$.
\par
\medskip
{\bf Corollary.} {\it For a difference tensor } $K${\it of
a Blaschke structure} $(\nabla,h,S)$ {\it the following equation holds:}

$$ [K_{X},K_{Y}]Z = -J(h(Y,Z)X - h(X,Z)Y)$$
\par
\medskip
{\bf Proof:} It follows from lemmas A and A1.$\diamondsuit$
\par
\medskip
{\it Remark.} It is clear that $\{ x:\epsilon (x) = \alpha  \} \cap  M$ is an open set  for
$\alpha  \in  \{- 1,1\}$. Let us note also that   relations  (4.4)  could  be
obtained from (4.3) by replacing $h$ by $- h$, which corresponds to  a
change of orientation of M. Hence we can restrict ourselves to the
investigation of surfaces $(M,f)$ satisfying (4.3) or (4.5).
\par
\bigskip
{\bf Lemma  B1.}  {\it Let} $(M,f)$ {\it be an affine nonconvex  surface  with  induced
Blaschke structure and} $\{ E_{1},E_{2}\}$ {\it be a local    frame  on}
$M${\it . If}  $\{ E_{1},E_{2}\}$ {\it satisfies  relations} (4.3) {\it then}
\par
\medskip
$$ \align
\nh K(X,E_{1},E_{1})& = (X\lambda )E_{1} + 3\lambda  \omega (X)E_{2},\\
 \nh K(X,E_{1},E_{2})& = - 3\lambda  \omega (X)E_{1} - (X\lambda )E_{2}
  \tag 4.3'\\
\nh K(X,E_{2},E_{2})& = (X\lambda )E_{1} + 3\lambda  \omega (X)E_{2}
 \endalign $$
\noindent {\it with} $\omega  = \omega ^{1}_{2} = \omega ^{2}_{1}$
{\it the connection form defined by }$\nh_{X}E_{i} = \omega ^{j}_{i}(X)E_{j}$.
{\it If}
$\{ E_{1},E_{2}\}$ {\it satisfies  relations} (4.5) {\it then the equations below hold}
\par
$$
\gather
\nh K(X,E_{1},E_{1}) = 0, \nh K(X,E_{1},E_{2}) = 0, \tag 4.5'\\
 \nh K(X,E_{2},E_{2}) = ((Xa) - 3a \omega (X))E_{1}.
\endgather $$
\noindent {\it with}
$\omega  = \omega ^{1}_{1} = - \omega ^{2}_{2}$
\par
\medskip
{\bf Corollary  1.}  {\it Let} $(M,f)$ {\it be an affine sphere and}
 $U = \{ h(C,C) \neq
0\}$. {\it Then for every point} $x_{0} \in  U$ {\it there  exists  a  chart} $(W,x_{1},x_{2})$
{\it with an associated local frame} $\{\partial _{1},\partial _{2}\}
 ( \partial _{i} = \frac{\partial}{\partial {x}_{i}})$ {\it such that}
  $x_{0} \in
W, \alpha $ = sgn $h(C,C)$ {\it is constant on} $W$
{\it and the equations below hold:}
\par
$$
E_{1} = e^{u}\partial _{1} , E_{2} = e^{u}\partial _{2} \tag E2
$$
$$
 h = \alpha  e^{-2u}(( dx_{1})^{2} - (dx_{2})^{2}), \tag h'
$$
$$
 K(\partial _{1},\partial _{1}) = e^{2u}\partial _{1}, K(\partial _{2},\partial _{2}) = e^{2u}\partial _{1},\tag K1
$$
$$
K(\partial _{1},\partial _{2}) = - e^{2u}\partial _{2}
$$
\noindent {\it and}
\par
$$ \gather
\nabla _{\partial _{1}}\partial _{1} = (- \partial _{1}u + e^{2u} )\partial _{1} -
\partial _{2}u \partial _{2}\\
 \nabla _{\partial _{1}}\partial _{2} = - \partial _{2}u \partial _{1}
 - (\partial _{1}u + e^{2u}) \partial _{2} \tag $\nabla 1$\\
\nabla _{\partial _{2}}\partial _{2} =
( - \partial _{1}u + e^{2u}) \partial _{1} - \partial _{2}u \partial _{2}
\endgather
$$
$$
 \Delta _{0} u = e ^{-2u} (\alpha  H + 2e^{6u}) \tag $\Lambda 1 $
$$
\noindent {\it or equivalently}
\par
$$
 \Delta _{0} \ln  \lambda  = 3 \lambda  ^{-\frac 23}
 ( \alpha  H + 2\lambda ^{2})
$$
\noindent {\it where}
$\Delta _{0} = \partial ^{2}_{1} - \partial ^{2}_{2}$
{\it is the standard Laplacian on} $( R^{2},dx^{2}_{1} - dx^{2}_{2})$
\par
\noindent $\lambda  = \frac14\sqrt{\alpha h(C,C)}$ {\it and}
$u = \frac13\ln \lambda $.
\par
\medskip
{\bf Proof:}  We omit the proof as it just the same as the  proof  of
Corollary below Lemma B.$\diamondsuit$
\par
\medskip
{\bf Corollary  2.}
 {\it Let us assume that} $(M,f)$  {\it is  an  affine  sphere,
with  indefinite affine metric} $h$ {\it and affine  mean  curvature} $H${\it ,
for which} $h(C,C) = 0$ {\it on} $M${\it . Then the Gauss curvature} $K_{h}$ {\it of} $(M,h)$ {\it is
constant and equals} $H${\it . Let us define} $U$ {\it =} int $\{ x \in  M: K_{x} = 0 \}$ {\it and}
$V = \{ x \in  M:K_{x} \neq  0 \}$ {\it where} $K$ {\it is the difference tensor  of} $(M,f)${\it .
An affine surface} $(U,f)$ {\it is a quadric and  for  every  point  of} $V$
{\it there exist a local frame} $\{E_{1},E_{2}\}$ {\it and  local coordinates} $(W,x_{1},x_{2})$
{\it such that the following equations hold:}
\par
$$
\gather
 E_{1} = e^{u} \partial _{1}\qquad E_{2} = \partial _{2} \tag E2\\
K(\partial _{1},\partial _{1}) = 0, K(\partial _{1},\partial _{2}) = 0,
 K(\partial _{2},\partial _{2}) = e^{u} \partial _{1}  \tag K2\\
h(\partial _{1},\partial _{1}) = h(\partial _{2},\partial _{2}) = 0 ,
 h(\partial _{1},\partial _{2}) = e^{-u}
\endgather
$$
\par
$$
 \partial _{1}(\partial _{2}u) = H \exp (- u) \tag $\Lambda 2$
$$
\noindent {\it In the introduced coordinates the immersion} $f$ {\it in a case} $H \neq  0$  {\it has
the form:}
\par
$$
\qquad f(z_{1},z_{2}) = \frac1H (\partial _{2}u(z_{1},z_{2})\xi (z_{2}) - \xi ^{'}(z_{2}))
\tag f$$
\noindent {\it where}
$\xi  : \Bbb{R} \rightarrow  \Bbb{R}^{3}$ {\it is a  smooth  curve  in}
$\Bbb{R}^{3}$  {\it defined  by  the
differential equation}
\par
$$
 \xi ^{'''} = a \xi ' + b \xi      \tag $\xi$
$$
\noindent {\it with the initial conditon} $\det (\xi ,\xi ',\xi '')_{0} = H$, {\it where} $a,b$
 {\it are some smooth functions such that} $a' - 2b = 2H.$
\par
\medskip {\it In a case} $H = 0 \ (V,f)$ {\it is equiaffinely equivalent to the  graph  of
a function} $z(x,y) = xy + \phi (y)$ {\it where} $\phi $ is {\it any smooth function such
that} $\phi ''' \neq  0.$
\par
\medskip
{\bf Proof:}  We consider only the case $H \neq  0$ as the case $H = 0$  is
well known (see \cite{M-R}). From lemma A1  there  exists  local  frame
$\{\bar{E}_{1},\bar{E}_{2}\}$ such that  equations (4.5) hold
 for  some  function $a \in
C^{\infty }(M)$. Let us denote by $V$ the set $\{ x:a(x) \neq  0 \}$.
 Let us define on
$V$ a local frame $\{ E_{1},E_{2} \}$ as follows
$E_{1} = a^{-\frac13} \bar{E}_{1}, E_{2} = a^{\frac13} \bar{E}_{2}$.
Equations (4.5) hold for $\{ E_{1},E_{2} \}$ with $a = 1$. We also have:
$$
\gather
\nh K(X,E_{1},E_{1}) = 0,\nh K(X,E_{1},E_{2}) = 0,\tag 4.6\\
 \nh K(X,E_{2},E_{2}) = - 3\omega (X) E_{1}
 \endgather
$$
\noindent where  $\nh_{X}E_{1} = \omega (X)E_{1}$, $\nh_{X}E_{2} = - \omega (X)E_{2}$.
From (4.6) it follows that   $\nh K$
is symmetric if and only if $\omega (E_{1}) = 0.$
Let us denote $\alpha  := \omega (E_{2})$. The equations
\par
$$
\nh _{E_{1}}E_{1} =\nh_{E_{1}}E_{2} = 0 ,\nh_{E_{2}}E_{1} = \alpha  E_{1},
 \nh_{E_{2}}E_{2} = - \alpha  E_{2} \tag 4.7
$$
\noindent hold on V. Let us introduce coordinates $y_{1},y_{2}$ such that
\par
$$
 E_{1} = \phi  \partial _{1} , E_{2} = \psi  \partial _{2} \tag 4.8
$$
\noindent From (4.7) it is clear  that $\partial _{1}\psi  = 0.$  Let  us  change  coordinates
taking $z_{1} = y_{1}, z_{2} = \Psi (y_{2})$ where $\Psi  = \int  \psi ^{-1}$.
 In new coordinates equations (4.8) hold with $\psi  = 1$
and we also have:
$$
\nh_{\partial _{1}}\partial _{1} = - \partial _{1}
\ln\mid \phi\mid \partial _{1},
\nh_{\partial _{1}}\partial _{2} = 0,
\nh_{\partial _{2}}\partial _{2} =
 - \partial _{2}\ln\mid\phi\mid\partial _{2}$$
and $\alpha  = \partial _{2}\ln \mid \phi \mid$. Let us define
$u := \ln \mid\phi\mid $. Notice that $[ E_{1},E_{2}] =
-\alpha E_{1}$. Hence $d\omega (E_{1},E_{2}) = E_{1}\omega (E_{2}) =
\phi  \partial _{1}\alpha  = \phi  \partial _{12}\ln \mid \phi \mid$ .
On the other
hand $d\omega (E_{1},E_{2}) = K_{h} = H$.
We get the equation $\partial _{1}(\partial _{2}u)$  = $H$exp$(-u)$.
Let us assume $H \neq  0.$ We can assume without loss of generality that
$- Hf$ is an affine normal  field  for $(M,f)$.  In  the  coordinates
$(W,z_{1},z_{2})$ the immersion $f$  satisfies  the  following  fundamental
system of equations:
$$
\align
\partial ^{2}_{1}f &= - \partial _{1}u \partial _{1}f\\
\partial ^{2}_{2}f &= e^{u} \partial _{1}f - \partial _{2}u \partial _{2}f
\tag F\\
\partial ^{2}_{12}f &= - e^{-u} H f
\endalign
$$
\noindent From the first of equations (F) it follows that $\partial _{1}f = e^{-u}\xi (z_{2})$ for
some smooth function $\xi : R \rightarrow   R^{3}$. From the third equation  we  get
$f(z_{1},z_{2}) = \frac1H ( \partial _{2}u \xi (z_{2}) - \xi ^{'}(z_{2}))$ and
\par
$$
\partial _{2}f = \frac1H ( \partial ^{2}_{2}u \xi  + \partial _{2}u \xi ^{'} - \xi  ^{''}).
$$
\noindent It is not difficult to see that the second equation  is  satisfied
iff $\xi $ satisfies the following equation:
\par
$$
\xi ^{'''} = ( 2 \partial ^{2}_{2}u + (\partial _{2}u)^{2})\xi ^{'}
+ ( \partial ^{3}_{2}u + (\partial _{2}u)(\partial ^{2}_{2}u) - H )\xi
$$
\noindent Let us note that in fact it is an ordinary  differential  equation
of the third order. If we denote by $a,b$ the functions $( 2 \partial ^{2}_{2}u +
(\partial _{2}u)^{2}),( \partial ^{3}_{2}u +
(\partial _{2}u)(\partial ^{2}_{2}u) - H )$ respectively then
$2b = \partial _{2}a - 2H$
and $a,b$ do not depend on $z_{1}$.  The  last  statement  follows  from
equation $(\Lambda 2)$. If a mapping $\xi $ satisfies equation $(\xi )$  then  it  is
obvious  that $\det (\xi ,\xi ^{'},\xi ^{''})$  is  constant.  One  can  check  that
$\det (\xi ,\xi ^{'},\xi ^{''}) = H$ as $f$ is a Blaschke immersion.
On the other  hand
if $u$ satisfies the Liouville equation
 $(\Lambda 2)$ on $\Omega  \subseteq   R^{2}$ and we choose
a solution $\xi $ of equation $(\xi )$ such that
 $\det (\xi ,\xi ^{'},\xi ^{''}) = H$  then  an
immersion $f$ given by  equation (f) is a Blaschke immersion  and
$(\Omega ,f)$ is an affine sphere with an affine normal $- Hf$.$\diamondsuit$
\par
\bigskip
{\it Remark.} If we change coordinates as follows:
\par
$$
x_{1} = \frac1H \partial _{2}u(z_{1},z_{2})\qquad x_{2} = z_{2}
$$
\noindent then it is clear that in new coordinates   equations $(f1)$
\par
$$
f(x_{1},x_{2}) = x_{1}\xi (x_{2}) - \frac1H \xi ^{'}(x_{2}) \tag f1
$$
\noindent and $(\xi )$ hold. It is easy to check that if we define an
immersion $f$
by equation $(f1)$ where $\xi $ satisfies equation $(\xi )$
together with  the
initial condition and $a,b$ are {\it any} smooth  functions $a,b
\in  C^{\infty }(\Bbb{R})$
then $(\Bbb{R}^{2},f)$  is  an  affine  sphere  with  vanishing  Fubini-Pick
invariant and nonzero affine mean curvature $H$. $( \Bbb{R}^{2},f)$ is a quadric
iff $a' = 2b.$
\par
\medskip
{\bf Theorem   3.}  {\it Let a function} $u$ {\it satisfies  equation}
 $(\Lambda 1)$  {\it on  an
open simply connected set} $\Omega  \subseteq  \Bbb{R}^{2}${\it .
Then there exists  an  immersion}
$f_{u}: \Omega  \rightarrow   \Bbb{R}^{3}$ {\it such that}
$(\Omega ,f_{u})$ {\it is an affine  sphere  with  affine
mean  curvature} $H$  {\it and  indefinite  affine  metric   h.   The
Fubini-Pick invariant of} $(\Omega ,f_{u})$  {\it equals} $16\alpha e^{6u}$.
  {\it In  the  standard coordinates on} $\Omega $ {\it  equations}
   (h '),$(K1),(\nabla 1)$ {\it hold where} $h$ {\it is  an
affine metric,} $K$ {\it is the difference tensor and} $\nabla $
  {\it is the  Blaschke connection for} $(\Omega ,f_{u})${\it .}
\par
{\it Let} $(M,f)$ {\it be  an  affine  surface  with  induced  Blaschke
structure which is an affine sphere with an affine mean  curvature}
$H${\it . Then} $M = U \cup \bar{V}= \bar{U} \cup  V$ {\it where} $U,V$
{\it are open  subsets  of} $M$  {\it such
that} $U \cap  V =\emptyset $  {\it ,} $U \subseteq  \{ x:h(C,C) = 0 \}$ {\it and}
$V = \{ x:h(C,C) \neq  0 \}$.
{\it Around any point} $x_{0} \in  V$  {\it there exist}
{\it a  local  chart} $(W,x_{1},x_{2})$
{\it such that the  function}
$u = \frac13\ln \lambda \circ \gamma ^{-1}${\it ,  where}
$\lambda  = \frac14\sqrt{\alpha h(C,C)}$ {\it ,
satisfies on} $\gamma (W)$ {\it  equation} $(\Lambda 1)$ {\it where}
 $\gamma (p) = (x_{1}(p),x_{2}(p)), \alpha $
{\it =} sgn$h(C,C)$ {\it and} $(W,f)$ {\it is  equiaffinely  equivalent  to  an  affine
sphere in} $ R^{3}$ {\it given by} $(\gamma (W),f_{u})$. {\it We also have}
$U = U_{0} \cup  U_{1}$,  {\it where}
$U_{0}$ = int $\{ x:K_{x} = 0 \}$ {\it and}
$U_{1} = \{ x:K_{x} \neq  0 \}\cap \{ x:h(C,C) = 0\}$.
{\it A surface} $(U_{1},f)$ {\it is locally equiaffinely equivalent to  one  of  the
surfaces described in the Corollary} 2 {\it and a surface} $(U_{0},f)$  {\it is  a
nonconvex quadric.}
\par
\bigskip
{\bf Proof:}  We omit the proof as it is  just  the  same  as  the
proof of the Theorem 2.
\par
\medskip
{\it Remark.}  Let us note that the equation ($\Lambda1$)
is equivalent to the equation
$$\Delta_0\Psi=e^{2\Psi}+\e e^{-\Psi} \tag $\Delta_5'$ $$
where $\e=\text{sgn}\a H\in\{-1,0,1\}$ and
$\Psi(x_1,x_2)=2u(\frac{x_1}a,\frac{x_2}a)-\ln b$ where
$a=(4\mid H\mid)^{\frac13}$, $b=\frac{\mid H\mid}2)^{\frac13}$.
if $H\ne 0$ and $\Psi=2u(\frac{x_1}2,\frac{x_2}2)$ if $H=0$.

\par
\medskip
{\bf Corollary  3.}  {\it Let} $(M,f)$ {\it be an affine sphere in}
 $\Bbb{R}^{3}$ {\it with  an  affine
hyperbolic metric} $h$ {\it and  Fubini-Pick invariant} $h(C,C)$.
 {\it Let  us
assume that} $h(C,C)$ {\it is different from} 0  {\it on} $M$.
 {\it Then for every} $x_{0} \in
M$ {\it there exist a neighborhood} $V$ {\it of} $x_{0}$
{\it and a one parameter family of
affine immersions} $f_{a}, a \in  O(1,1)${\it , such that} $(V,f_{a})$
 {\it is  an  affine
sphere and each immersion} $f_{a}$ {\it has the same induced affine metric} $h$
{\it and the same Fubini-Pick invariant} $h(C,C)${\it .  Every  affine
  sphere
immersion whose  induced Blaschke structure has the same as }$(M,f)$
{\it  affine metric} $h$
{\it and  Fubini-Pick  invariant} $h(C,C)$
 {\it is  locally  equiaffinely
equivalent to one of immersions} $f_{a}$.
\par
\medskip
{\bf Proof}: The symmetry group of  equation $(\Lambda 1)$ is  an  affine
orthogonal group $AO(1,1)$ of  all  affine  transformations  of  the
form:
\par
$$
 g = \pmatrix
 \cosh \alpha &\epsilon \sinh \alpha &a \\
\sinh\alpha &\epsilon\cosh \alpha &b \\
 0&0&1
 \endpmatrix \in  AO(1,1), \alpha ,a,b \in \Bbb{R} \tag g1
$$
\noindent Let us define  in the coordinates  described  in  the  Theorem 3
a difference tensor $K_{\alpha }$ by the formula:
\par
$$
\gather
K_{\alpha }(\partial _{1},\partial _{1}) = e^{2u}( \sinh  3\alpha
 \partial _{1} - \epsilon  \cosh  3\alpha  \partial _{2})\\
 K_{\alpha }(\partial _{2},\partial _{2}) = e^{2u}
 ( \sinh  3\alpha \partial _{1} - \epsilon  \cosh  3\alpha  \partial _{2}) \tag $\alpha 1 $\\
K_{\alpha }(\partial _{1},\partial _{2}) = e^{2u}
(\epsilon  \cosh  3\alpha  \partial _{1} - \sinh  3\alpha  \partial _{2})
\endgather $$
\par
\noindent Let $\nh$  be the Levi-Civita connection of  an  metric $h$
 defined  by
(h) and define a connection $\nabla _{\alpha }$ as $\nabla _{\alpha }$ =
$\nh  + K_{\alpha }$ where $K_{\alpha }$ is  given
by the formulas $(\alpha 1)$. Then as before one can prove  the  existence
of an immersion $f_{a}$, where
\par
$$
a = \pmatrix
\cosh\alpha &\epsilon\sinh\alpha \\
\sinh\alpha &\epsilon \cosh\alpha
\endpmatrix
$$
\noindent with a Blaschke connection $\nabla _{\alpha }$ and an affine metric
 h. The rest  of
the proof is the same as the  proof  of  Corollary  following  the
Theorem 2.$\diamondsuit$
\par
\medskip
{\bf Proposition 1.} {\it Let a function} $u$ {\it satisfies on} $\Omega $
{\it  equation} $(\Lambda 1)$
{\it and let} $(\Omega ,f_{u})$ {\it be an affine sphere with induced  affine metric
  and induced  Blaschke  connection
given by formulas} $(h')$ {\it and} $(\nabla 1)$. {\it Let} $(\Omega ,f_{a})$
 {\it be as above,} $f_{id} = f_{u}${\it .
Then} $(\Omega ,f_{a})$ {\it is affinely equivalent to} $(\Omega ,f_{u})$
{\it iff} $u = u\circ g$, {\it where}  a
{\it is given by the formula} (a) {\it and} $g$ {\it is given by} (g)
{\it for some} $a,b \in \Bbb{R}${\it .}
\par
\medskip
{\bf Proof:}  The proof is the same as  the  proof  of
the Proposition.$\diamondsuit$
\par
\medskip
{\it Remark.}  Let us  note  that  if $h(C,C)$  is  constant  and
different from 0 then all hypersurfaces $( \Omega ,f_{a})$  are  equiaffinely
equivalent. In particular an  affine sphere with constant
Fubini-Pick invariant $h(C,C) \neq  0$ and indefinite affine metric  is  characterized  uniquely  by
$h(C,C)$. Every non-convex affine  sphere  with  constant  nonzero  Fubini-Pick
invariant is equiaffinely  equivalent  to  an  open  part  of  the
surface $Z(X^{2}+Y^{2}) = \pm  \frac1c$ where $c \neq  0$ depends  only  on $h(C,C)$. The
fundamental system of equations for $f$ is very simple in  the  case
of constant nonzero $h(C,C)$ and it is easy to see that (in a  chart
for which $\partial _{i} = E_{i})$

$$
f(x_{1},x_{2}) = ( e^{-\lambda x_{2}} \cos ( \sqrt 3 \lambda x_{1}),
- e^{-\lambda x_{2}} \sin (\sqrt 3\lambda x_{1}),\pm \frac1c
e^{2\lambda x_{2}} )
$$
where $c = \frac{3\sqrt 3}{128} h(C,C)^{2}$ (see also \cite{M-N}).  The  surface
$M = \{(X,Y,Z) :Z(X^{2}+Y^{2}) = \pm  \frac1c \}$ is affinely homogeneous.
 It is  easy  to
check that $( \Bbb{R}^{2},f )$ is an orbit of the point $(1,0,\pm  \frac1c)$ by
 the group
$G =  \Bbb{R}\oplus  S^{1}$ of equiaffine transformations
\par
$$
\pmatrix
e^{-a} \cos b & e^{-a} \sin b& 0 \\
-e^{-a} \sin b & e^{-a} \cos b & 0 \\
0 & 0 & e^{2a}
\endpmatrix (a,b) \in \Bbb{R} \oplus \frac{\Bbb{R}}{2\pi \Bbb{Z}}
$$
\par
\medskip
\centerline{{\bf References.}}
\par
\medskip
\cite{B-C} `Solitons' edited by R.K. Bullough and P.J. Caudrey, Springer
Verlag Berlin Heidelberg New York 1980
\par
\medskip
\cite{D}  Dillen,F.:`{\it Equivalence   theorems   in   affine
differential geometry.}' Ge\-om. Dedi\-ca\-ta {\bf 32}, (1989), 81 - 92.
\par
\medskip
\cite{D-N-V} F.Dillen, K.Nomizu, L.Vrancken `{\it Conjugate Connections  and
Radon's Theorem in  Affine Differential  Geometry}'   Mh.Math.{\bf 109},
221 -235,(1990).
\par
\medskip
\cite{J-1} W.Jelonek.: `{\it Affine locally symmetric surfaces}' Geom.
Dedicata {\bf 44}, 189-221, (1992).
\par
\medskip
\cite{J-2} W.Jelonek.:`{\it The local structure of affine spheres in} $ R^{3}$' Preprint 505,
Polish Academy of Sciencies, (1992).
\par
\medskip
\cite{M-N} Magid, M., Nomizu, K.`{\it On Affine Surfaces whose  Cubic  Forms
are Parallel Relative to the  Affine  Metric.}Proc.  Japan.  Acad.,
{\bf 65},Ser.{\bf A} (1989),215-218.
\par
\medskip
\cite{M-R} Magid, M.,Rayan, D.`{\it Flat affine spheres in} $ R^{3}.'$Geom.
Dedicata
{\bf 33}, No.3, (1990), 277-288.
\par
\medskip
\cite{N} V.V.Nesterenko `{\it On the geometric origin of the equation}
$\phi,_{11}-\phi,_{22}=e^{\phi}-e^{-2\phi}$' Letters in Mathematical
Physics {\bf 4}, 451-456, (1980)
\par
\medskip
\cite{N-P-1} K.Nomizu, U. Pinkall.: `{\it On a certain  class  of  homogeneous
projectively flat manifolds}' T\^ohoku Math.Journ. {\bf 39}, 407-427,(1987)
\par
\medskip
\cite{N-P-2} K.Nomizu, U. Pinkall.:`{\it On   the   Geometry    of    Affine    Immersions}'
Math. Z. {\bf 195}, 165-178, (1987).
\par
\medskip
\cite{N-P-3} Nomizu,K.,Pinkall,U. `{\it Cayley Surfaces in Affine Differential
Geometry.}' T\^ohoku Math.J. {\bf 41},(1989), 589-596.
\par
\medskip
\cite{P} Podesta,F.:`{\it Projectively flat affine surfaces in} $\Bbb{A}^3$' Proc.
Amer. Math. Soc. {\bf 119}, 255-260 (1993).
\par
\medskip
\cite{S-W} Simon,U.,Wang,Ch.,  `{\it Local  Theory  of  Affine
2-Spheres.}' Proc. of Symposia in Pure Math., part 3, {\bf 94}
585-598, (1993).
\par
\medskip
\cite{S} /Slebodzi/nski, W.: {\it `Sur quelques problems de la theorie
des surfaces de l'espace affine'} Prace Mat. Fiz. {\bf 46},291-345,
(1939).

\par
\medskip
\noindent Institute of Mathematics
\par
\noindent Technical University of Cracow
\par
\noindent Warszawska 24
\par
\noindent 31-155 Krak/ow,POLAND.

\end